% This file was written in AMS-TeX version 2.x
\input amstex

\magnification=\magstep1
\pageheight{22 truecm}
\pagewidth{16.5truecm}
%\hcorrection{-0.5 cm}
%\vcorrection{-0.5 cm}
\documentstyle{amsppt}
\loadbold
\NoBlackBoxes

\def\openN{I\negthickspace N}
\def\ldoublebracket{[\![}
\def\rdoublebracket{]\!]}
\def\edim{\operatorname{edim}}
\def\Ker{\operatorname{Ker}\,}
\def\Coker{\operatorname{Coker}\,}
\topmatter
\title
Report on the Torsion of the Differential Module
of an algebraic Curve
\endtitle
\rightheadtext{Torsion of the Differential Module}
\author
 Robert W. Berger
\endauthor
\address
{Fachbereich Mathematik, Universit\"at des Saarlandes, Postfach
151150, D 66041 Saarbr\"ucken, Germany}
\endaddress
\dedicatory
Dedicated to Professor Shreeram S. Abhyankar on the occasion of his
sixtieth
birthday
\enddedicatory
\abstract
There is a conjecture, that the torsionfreeness of the module
of differentials in a point of an algebraic or algebroid curve
should imply that the curve is non singular at that point. A
report on the main results is given.
\endabstract
\endtopmatter

\document
Let $k$ be a perfect field and  $R$  the local ring of a closed point
of an algebraic or algebroid curve over $k$. There is a conjecture
that  $R$  is regular if and only if the (universally finite)
differential module $\Omega_{R/k}$ is torsionfree. The nontrivial part
is
of course to show that for a singular point the torsion submodule
$\tau(\Omega_{R/k})$ of $\Omega_{R/k}$ is not zero. Although a solution
for the
general case is not in sight there are many special cases which
have been treated successfully. In all of these the conjecture has
been found to be true. It is the purpose of this paper to give a
survey on some of these results with hints concerning the proofs.
For simplicity let us assume for the following that $R$ is a
reduced complete analytic $k$-algebra of dimension one with maximal
ideal $\frak m$ and embedding dimension $n$, which then can be
represented in the form $R = k\ldoublebracket X_1,\dots
,X_n\rdoublebracket /I
=
 k\ldoublebracket x_1,\dots ,x_n\rdoublebracket$ ,
where $I$ is a reduced ideal in the formal power series ring
$k\ldoublebracket
X_1,\dots ,X_n\rdoublebracket$ .
We will also restrict ourselves to the case $char\, k = 0$ ,
although many of the results are also valid for perfect ground
fields.
\newline
One can distinguish several cases:
\head 1.
Conditions on the number of generators of $I$ .
\endhead
Let $d(R) = \mu(I)-(n-1)$ denote the {\it deviation\/} of $R$ , where
$\mu(I)$ denotes the minimal number of generators of $I$ . $R$ is called
a {\it complete intersection\/} if $d = 0$ and an {\it almost complete
intersection\/} if $d\le 1$ . In [Be2] the cases $d\le 1$ were solved if
$R$ is a domain. This was generalized to $d\le 3$ in the reduced case
by Ulrich [Ul1], [Ul2].
\par
Denote by $S$ the integral closure of $R$ in its full ring of
quotients $K$ , and let  $D:S\rightarrow SDS = \Omega_{S/k}$  and
$d:R\rightarrow RdR = \Omega_{R/k}$  be the universally finite
derivations of $S$ and $R$  over $k$ respectively.
Since $RdR$ and $SDS$ are both of rank~1 and $SDS$ is torsionfree
(even free), the kernel of the canonical homomorphism
$\varphi : RdR \rightarrow SDS$ is $T := \tau(\Omega_{R/k})$, so that
we have an exact sequence
$$0\rightarrow T\rightarrow RdR\rightarrow SDS\rightarrow
SDS/RDR\rightarrow
0\,.$$
Let $x\in R$  be a normalizing parameter and $ s = k\ldoublebracket
x\rdoublebracket$ . Then $R$
and $S$ are finite $s$--modules and the universally finite derivations
$\delta\, : R\rightarrow R\delta R = \Omega_{R/s}$ and
$\Delta\, : S\rightarrow S\Delta S = \Omega_{S/s}$ coincide
with the universal derivations of $R$ and $S$ over $s$
respectively. $R\delta R$ and $S\Delta S$ are finitely generated torsion
modules and therefore have finite lengths as $s$--modules. Further we
have $R\delta R = RdR/Rdx$ and $Rdx\cap T=0$ and therefore $SDx/RDx\cong
S/R$ .
Let $\alpha$ be the natural
map $R\delta R\rightarrow S\Delta S$ induced
by the inclusion $R\hookrightarrow S$ . One gets a commutative diagram
with
exact rows and columns:
$$
\CD
{}  @.    {}    @.     0      @.        0               @.    {} \\
@.        @.           @VVV             @VVV                  @. \\
{}  @.    0     @.     T      @.        {\Ker}\,\alpha     @.    {} \\
@.        @VVV         @VVV             @VVV                  @. \\
0   @>>>  Rdx   @>>>   RdR    @>>>      R\delta R       @>>>  0  \\
@.        @VVV         @V{\varphi}VV    @V{\alpha}VV          @. \\
0   @>>>  SDx   @>>>   SDS    @>>>      S\Delta S       @>>>  0  \\
@.       @VVV         @VVV             @VVV                   @. \\
{}  @.    S/R   @.   SDS/RDR  @.      {\Coker}\,\alpha     @.    {} \\
@.        @VVV         @VVV             @VVV                  @. \\
{}  @.    0     @.     0      @.        0               @.    {} \\
\endCD
$$
By the snake lemma one obtains an exact sequence
$$0\rightarrow T\rightarrow Ker\,\alpha\rightarrow S/R\rightarrow
SDS/RDR\rightarrow {\Coker}\,\alpha\rightarrow 0$$
and from this for the lengths as s-modules:
$$l(T)=l(Ker\,\alpha )-l({\Coker}\,\alpha )+l(SDS/RDR)-l(S/R)\,.$$
But $l(Ker\,\alpha )-l({\Coker}\,\alpha )=l(R\delta R)-l(S\Delta S)$ ,
so that we have
\newline
\proclaim{Lemma 1}
$$l(T)=l(R\delta R)-l(S\Delta S)+l(SDS/RDR)-l(S/R)\,.$$
\endproclaim
\subhead 1.1.
The case $\boldkey d \boldkey ( \boldkey R \boldkey )
\boldsymbol\leq \boldkey 1$\
\endsubhead
One can represent $R\delta R=F/U$ where $F$ is a free $R$--module and
$U$
is generated by at most $rank(F)+1$ elements. Therefore by \cite{Be2,
Satz 2}
 one has $l(R\delta R)=l(D(R\delta R)^{-1}/R)$\,, where $D(R\delta R)$
denotes the $0^{th}$ Fitting ideal of $R\delta R$, which is by
definition
also the $0^{th}$ K\"ahler Different ${\frak D}_K(R/s)$ of $R$ over
$s\,.$
By the general theory of the differents one has the following inclusions
for the Dedekind, K\"ahler, and Noether differents ${\frak D}_D$,
${\frak D}_K$, ${\frak D}_N$  and the complementary modules $R^*=\frak
C(R/s)$
and  $S^*=\frak C(S/s)\,:$
$$R\subseteq S\subseteq S^*\subseteq R^*\subseteq {\frak D}_D(R/s)^{-1}
={\frak D}_N(R/s)^{-1}\subseteq {\frak D}_K(R/s)^{-1}\, .$$
Therefore we have
$$
\align
l(R\delta R) =l({\frak D}_K(R/s)^{-1}/R) &\ge l(R^* /R) \\
             & =l(S^*/S)+l(R^*/S^*)+l(S/R) \\
             &=l(S^*/S)+2\cdot l(S/R) \text{\ by duality.}
\endalign
$$
Now $S$ is a principal ideal ring and $S\Delta S$  is of projective
dimension $1$ , so
\newline
$l(S\Delta S)=l(S/{\frak D}_K (S/s))=l(S/{\frak D}_D(S/s))=
l(S^*/S)\,,$ and therefore
\newline
$l(R\delta R)-l(S\Delta S)\ge 2\cdot
l(S/R)\,.$
\par
If $R$ is a complete intersection one even has $R^*={\frak D}_D
(R/s)^{-1}$ and
\newline
${\frak D}_K (R/s)={\frak D}_D(R/s)$  by \cite{Ku2, \S 10 and
appendix G,G1 and G2}.
\newline
Therefore $l(R\delta R)-l(S\Delta S)=2\cdot l(S/R)$  and we obtain
\proclaim{Theorem 1}  (\cite{Be2, Satz 7}): If $d\le 1$  then
$$l(T)\ge l(SDS/RDR)+l(S/R)$$
with equality if  $R$  is a complete intersection.
\par
Consequently:\ \  If $T=0$  then $R$ is regular.
\endproclaim
\par
There is another expression for $l(T)$ in this case due to Kunz
\cite{Ku1}:
\newline
Look at the natural homomorphism of $RdR$ into the $R$-module
$\omega_R$ of regular
\newline
differentials:
$$c_R\: RdR \rightarrow R^*\cdot Dx=\omega_R$$
induced by the inclusion $R \hookrightarrow Quot(R)$\,.
Obviously ${\Ker}\, c_R = T\, ,$ but also ${\Coker}\, c_R$ is
interesting in
this context: One has the inclusions
\newline
$R^*\cdot Dx\supseteq S^*\cdot Dx\supseteq RDR\supseteq RDx$
and so
$$
\align
l({\Coker}\, c_R)& =l(R^*\cdot Dx/RDR)=l(R^*/R)-l(RDR/RDx) \\
            & =l(R^*/S^*)
l(S^*/S)+l(S/R)+l(SDS/RDR)-l(SDS/SDx)-l(SDx/RDx) \\
            & =l(S/R)+l(S\Delta S)+l(S/R)+l(SDS/RDR)-l(S\Delta S)-l(S/R)
\\
            & =l(SDS/RDR)+l(S/R)\, .
\endalign
$$
If $R$ is a complete intersection then it follows with theorem 1
that $l(T) = l({\Coker}\, c_R)\, .$ If $R$ is only an almost complete
intersection, Kunz shows in  \cite{Ku1} that
\newline
$l(T) - l({\Coker}\, c_R) = l(\tau (I/I^2))\,$ ,  where $\tau (\ )$
denotes
the torsion submodule. For a complete intersection $\tau (I/I^2) = 0$ by
\cite{Ku1,theorem 1}, so  that we have
\proclaim{Theorem 1'} If $d\leq 1$ then
$$l(T) = l({\Coker}\, c_R) + l(\tau (I/I^2))\ .$$
\endproclaim
\subhead{1.2. The case $\boldkey d\boldkey (\boldkey R\boldkey )
\boldsymbol \leq\boldkey 3 $}
\endsubhead
In general it is not possible to give an inequality for the length
of an $R$-module using only the $0^{th}$ Fitting ideal. So the proof
from  1.1.  cannot be applied. (See the examples in  \cite{B2}.) But
Ulrich developed a more complex formula, using a chain of certain
determinantal ideals, by which a lower bound for the length of an
arbitrary finitely generated $R$-module can be computed. (\cite{Ul1,
Satz
2}, \cite{Ul2, Satz 1}).
\par
\noindent We start again with  Lemma 1 , but this time we choose the
parameter $x$ so that $SDR = SDx$. This is possible by  \cite{Ul1,
Hilfssatz 3}.  By Lemma 1 we have
$$
\align
l(T) &= l(R\delta R)-l(S\Delta S)+l(SDS/SDR)+l(SDR/RDR)-l(S/R) \\
     &= l(R\delta R)-l(S\Delta S)+l(\undersetbrace \ =S\Delta
S\to{SDS/SDx})+l(S\cdot RDR/RDR)-l(S/R)\text{\  by
 choice of $x$\, .}
\endalign
$$
Hence
\proclaim{Proposition 1} (\cite{Ul1, Satz 3})
$$
\align
l(T) &=l(R\delta R)-l(S/R)+l(S\cdot RDR/RDR) \\
     &\geq l(R\delta R)-l(S/R)\, .
\endalign
$$
\endproclaim
\definition{Definition 1} (\cite{Ul1, Def.3})
$R$  has {\it minimal torsion\/} if and only if
\newline
$l(T) = l(R\delta R) - l(S/R)\, .$
\enddefinition
\proclaim{Proposition 2} (\cite{Ul1, Bemerkung 3})  If $R$  has minimal
torsion
then
there exists a non zero divisor  $y\in {\frak m}$  with
$y\cdot {\frak m} = {\frak m}^2\, .$
If $R$ is singular then $T\neq 0\, .$
\endproclaim
\demo{Proof} If $R$ has minimal torsion then by proposition 1
$l(S\cdot RDR/RDR)=0\, ,$ hence $RDR = S\cdot RDR =
S\cdot Dx\cong S\,.$
$\Longrightarrow l_R({\frak m}/{\frak m}^2)=n\geq\mu_R (RDR)=\mu_R (S)=
l_R (S/{\frak m}\cdot S)\, .$  Now there is a non
zero divisor $y\in {\frak m}$ with  $S\cdot {\frak m} = S\cdot y$ since
$S$
is a principal ideal ring (\cite{Ul1, Hilfssatz 2}). Then
$l_R({\frak m}/y\cdot {\frak m}) = l_R(S/S\cdot y) = l_R(S/S\cdot {\frak
m})
\l_R({\frak m}/{\frak m}^2)\, $
hence ${\frak m}^2\supseteq y\cdot {\frak m}\supseteq {\frak m}^2\, .$
By theorem 8 we have then $T\neq 0\, ,$ if $R$ is singular.
\enddemo
Now by the above mentioned length formula of Ulrich, using also
the relations between the various differents and the fact that $R^*$
as in 1.1 is a canonical $R$-ideal, Ulrich shows in
\cite{Ul1, Satz 6}\,:
\proclaim{Proposition 3}  If $d\leq 3$ then
$l(R\delta R)\geq l(S/R) + l(S^*/S)\geq l(S/R)\, .$
\endproclaim
Now we can prove the torsion conjecture in this case:
\proclaim{Theorem 2}  (\cite{Ul1, Satz 6}, \cite{Ul2, Satz 3}):
If $d(R)\geq 3$ and $T=0$ then $R$ is regular.
\endproclaim
\demo{Proof}  By propositions 1 and 3  we have
$0=T\geq l(R\delta R) - l(S/R)\geq 0\, ,$ and hence $R$ has minimal
torsion.
Then by proposition 2  $T\neq 0$ if $R$ is not regular.
\enddemo
\remark{Remark 1} From propositions 1 and 3  we also get
$0 = l(T)\geq l(S^* /S)\geq 0\, ,$ and therefore $S^* = S$ . So $S$ over
$s$
is unramified.  If $k$ is algebraically closed and $R$ a domain,
it follows immediately that $S = s$ and therefore also $R = S$
without using theorem 7.
\endremark
\head{2. Exact differentials, maximal torsion and quasi homogeneous
singularities.}
\endhead
\subhead{2.1. Exact Differentials}
\endsubhead
A second class of curve singularities for which the conjecture is
true consists of all $R$ such that every differential of $RDR$ is
exact:
\proclaim{Theorem 3} (\cite{Po3, Theorem 3}):  Assume $RDR = DR\, .$
If $T = 0$ then $R$ is regular.
\endproclaim
For the proof we need the following lemma (\cite{G\"u2, Satz 1}), which
can be proved by a direct computation:
\proclaim{Lemma 2} Let  $A := k\ldoublebracket
Z_1,...,Z_q\rdoublebracket$ be a
formal power series
ring over $k$, $t\in \openN \, ,$  $A_t := A/(Z_1,\dots ,Z_q)^t\, .$
Then
${\dim}_k\Omega_{A_t/k} = (t-1)\cdot \binom{q+t-l}{t}\ .$
\endproclaim
\demo{Proof of theorem 3} Let $x\in {\frak m}\backslash {\frak m}^2$  be
a non
zero
divisor of $R$ , $\bar R := R/x$ and $\bar d:\bar R\rightarrow \bar
R\bar d\bar
R$
the universally finite derivation of $\bar R$ over $k$ .
Then ${\edim}\,\bar R = {\edim}\, R\,-1 = n-1\, ,$
so that we can represent
$\bar R = P/{\frak b}$\, ,\  $P := k\ldoublebracket
Y_1,...,Y_{n-1}\rdoublebracket$
 a formal power series ring with maximal ideal
 ${\frak M} = (Y_1,...,Y_{n-1})$ and $\frak b$ an ideal with
${\frak b}\subseteq {\frak M}^2\, .$
Denote by $\partial\,:\,P/{\frak M}^2 \rightarrow \Omega_{\frac{
P/{\frak
M}^2}{k}}$
the universally finite derivation of $P/{\frak M}^2$ over $k$.
We obtain the following commutative diagram with exact rows of canonical
maps:
$$
\CD
R                   @>>>   P/{\frak M}^2                   @>>>  0 \\
@V{\bar d}VV               @V{\partial}VV                       @.\\
\bar R\bar d\bar R  @>>>  \Omega_{\frac{P/{\frak M}^2}{k}}  @>>>  0 \\
\endCD
$$
This yields a surjection $\bar R\bar d\bar R/Im\,\bar d\rightarrow
\Omega_{\frac{P/{\frak M}^2}{k}}/Im\,\partial\rightarrow 0\, .$ Now by
assumption we have  $0 = T = {\Ker}(RdR\rightarrow RDR)$ and $RDR = DR$.
Hence $RdR = dR$ and consequently $\bar R\bar d\bar R = RdR/Rdx = \bar
d\bar R
 = Im\,\bar d\, .$ \ $\Longrightarrow$
$Im\,\partial = \Omega_{\frac{ P/{\frak M}^2}{k}}\, .$\ $
\Longrightarrow$
\ $n={\dim}_k\, P/{\frak M}^2 ={\dim}_k\, \Omega_{\frac{P/{\frak
M}^2}{k}}+
{\dim}_k\, Ker\, \partial\geq \binom{n}{2}+1$ by Lemma 2. $
\Longrightarrow
n\leq 2 \, .$ But
$n = 2$ is not possible, because then $R$ is a plane curve
singularity and therefore $T\neq 0$ by 1.1 .  Then $n = 1$ and
consequently $R$ is regular.
\enddemo
\par
When is the condition $RDR = DR$ satisfied?
\subhead{2.2.  Maximal Torsion}
\endsubhead
Consider the universally finite derivation $S\rightarrow SDS$ of $S$
over
$k$\, . Since $S$ is a direct product of formal power series rings
$S_i\cong k_i\ldoublebracket t_i\rdoublebracket\, ,$ where $k_i$ is an
algebraic extension of $k$ ,
$D$ is the direct product of the formal derivations
$D_i\: f(t_i)\mapsto f'(t_i)\cdot D_it_i\, .$ It follows by formal
integration that the $D_i$ are surjective and hence also $D$ is
surjective.
This induces a surjective $k$-linear map $\tilde D : S/R\rightarrow
SDS/RDR\,
,$
and so $l(S/R) = {\dim}_k(S/R)\geq {\dim}_k(SDS/RDR) = l(SDS/RDR)\, .$
Together with lemma 1  we get
$$
l(T)=l(R\delta R)-l(S\Delta S)-[\undersetbrace \ge 0
\to{l(S/R)-l(SDS/RDR)}]
\le l(R\delta R)-l(S\Delta S)\, .
$$
\definition{Definition 2} (\cite{Po3}):  R has {\it maximal torsion}  if
and
only if
$$l(T) = l(R\delta R) - l(S\Delta S)$$
or, equivalently, if and only if
$$l(S/R) = l(SDS/RDR)\,.$$
\enddefinition
If $R$  has maximal torsion the conjecture is true:
\proclaim{Theorem 4} (\cite{Po3, Theorem 1})
\newline
If  $R$  has maximal torsion then $RDR = DR$\, .
\newline
Consequently:\ \ \ If $T=0$  then  $R$  is regular.
\newline
If $k$ is algebraically closed and $R$ is a domain then also the
converse is
true:
\newline
\indent If  $RDR = DR$  then  $R$  has maximal torsion.
\endproclaim
\demo{Proof}  As shown above the $k$-linear map $D\: S\rightarrow SDS$
is
surjective and therefore
$$
\align
l(S/R) &= {\dim}_k(S/R)\ge {\dim}_k(S/(R + {\Ker}\, D)) =
{\dim}_k(SDS/DR) \\
       &\ge {\dim}_k(SDS/RDR) = l(SDS/RDR)\, .
\endalign
$$
If $R$ has maximal torsion then  $l(S/R)=l(SDS/RDR)$  and therefore $DR
=
RDR\,.$
\newline
Assume now that $k$ is algebraically closed and $R$ is a domain.
Then  ${\Ker}\, D = k \subseteq R$  and therefore
$l(S/R)=l(SDS/DR)=l(SDS/RDR)$
by
hypothesis. So $R$ has maximal torsion.
\enddemo
\remark{Remark 2} Let $k$  be algebraically closed and $R$ a
domain.
\newline
a) Let $ x \in \frak m\backslash{\frak m}^2$  be a superficial element
of $R$ .
Then one can find a uniformizing parameter $t$ for $S$ such that
$x=t^{m(R)}\,
,$ where  $m(R)$ denotes the multiplicity of $R$.  With $s\:=
k\ldoublebracket
x\rdoublebracket$  and
$S=k\ldoublebracket t\rdoublebracket$  we get $l(S\Delta
S)=l(SDt/SDt^{m(R)})=m(R)-1.$  Therefore $R$
has maximal torsion if and only if  $l(T)=l(R\delta R)-m(R)+1\, . $
This shows that definition 2  is equivalent to Ulrich's Definition 2
in \cite{Ul1}.
\newline
b) Zariski \cite{Za1} considers the case of an irreducible plane
algebroid
curve over $k$ and shows that $l(T)\le 2\cdot l(S/R)\,.$  Equality holds
if
and only if  $RDR=DR\,,$ which by theorem 4 means that $R$ has maximal
torsion.
He proves that this is the case if and only if the curve can be
represented
by a quasi homogeneous equation of the form $ Y^p - X^q =0$  with
$(p,q)=1\,.$
\endremark
For upper bounds of the torsion of a plane curve in terms of the
characteristic pairs see Azevedo \cite{Az, propositions 3 and 4}.
\subhead{2.3. Quasi homogeneous Singularities}
\endsubhead
\definition{Definition 3} (\cite{Sch, 9.8} ): $R$ is called {\it quasi
homogeneous}
if there exists a surjective $R$-module homomorphism
$\Omega_{R/k}\rightarrow\frak m$\,.
\newline
Let $\gamma =(\gamma_1,\dots,\gamma_n)\in {\openN}_0^n\,$ . A polynomial
$F =\sum\alpha_{i_1\dots i_n}\cdot X_1^{i_1}\dots X_n^{i_n} $
with coefficients in $k$ is called  quasi homogeneous of type $\gamma\,$
,if there is a $ d\in \openN$  such that  for all $ \alpha_{i_1\dots i_n}
\ne 0$
one has $ \sum\gamma_j\cdot i_j =d\,$. In this case  $d$ is called the
degree of~$F$\,.
\enddefinition
{}From \cite{K-R, Satz 2.1. and Satz 3.1.} one obtains:
\remark{Remark 3}: If $I$ is generated by polynomials, then  $R$ is
quasi homogeneous if and only if $I$ is generated by quasi homogeneous
polynomials of a fixed type $\gamma\,.$  For an irreducible $R$ this is
equivalent to $R$ being isomorphic to the analytic semigroup ring
$k\ldoublebracket H\rdoublebracket$ for the value semigroup
$H$ of $R$\, .
\endremark
\proclaim
{Theorem 5} (\cite{Po3, theorem 2}):
If $R$ is quasi homogeneous and $I$ generated by polynomials then $R$
has maximal torsion.
\par
Consequently:  If $ T = 0$  then $R$ is regular.
\endproclaim
\demo{Proof}: (For a different proof of the second assertion without the
assumption of $I$ being generated by polynomials see  Scheja
\cite{Sch}, Satz 9.8.) By hypothesis there is a $R$-linear map
$ \psi: RdR \rightarrow \frak m $ , which by  \cite{K-R}, proof of Satz
2.1,
can
be chosen such that $ \psi (dx_i) = \gamma_i x_i $ with $ \gamma_i \in
\Bbb N $
 for $ i = 1,...,n $. Since $RdR$  and $\frak m$  are both of rank 1
and
$\frak m $
is torsionfree we have $ {\Ker} \psi = T$.  We may assume that $ x =
x_1$  is
a superficial element of degree $1$ of $R$ .  Since $ Rdx \cap T=0$  we
have an exact sequence $ 0 \rightarrow T \rightarrow RdR/Rdx
\rightarrow \frak
m /R  x \rightarrow  0 $.
Therefore  $l(T) = l(RdR/Rdx) - l( \frak m /R x) = l(R \delta R) - m(R)
+ 1 $.
Now by remark 2 a) $R$ has maximal torsion.
\enddemo
\subhead 2.4.  The value semigroup
\endsubhead
Let, as before,  $R$  be a domain and  $k$  algebraically closed. Then
$S = k\ldoublebracket t\rdoublebracket$ is a discrete valuation ring.
Let
$\nu$  denote the normed
valuation with $\nu(t) = 1\,$. $ H =\{\nu (y)| y\in R\backslash 0\}$ is
called
the value semigroup of $R$\, .  Since $ SDS = SDt$  every $\omega\in
RDR$
is of the form  $\omega = z\cdot Dt$  with  $z\in S$,  and so we can
define
$\nu (\omega)\:=\nu (z)+1\,$.  This definition is independent of the
choice of
$t$\, . For all $y\in\frak m$  we then have $\nu (Dy) =\nu (y)$  and
therefore
$\nu (RDR)\supseteq\nu (\frak m)\,$,  but there may be elements
$\omega\in
RDR$
with $\nu (\omega) \notin\nu (\frak m)\,$ .  Yoshino calls them
exceptional
differentials (\cite{Yo, Def. 2.3.}). If $RDR = DR$  then obviously
$\nu (RDR) =\nu (\frak m)\,$. The converse follows from \cite{Za1, proof
of corollary 3}.
So we obtain
\proclaim{Theorem 6} (\cite{Ul1, Satz 4}):  Let $k$  be algebraically
closed,
and $R$  a domain. Then $\nu (RDR) =\nu (\frak m)$  if and only if  $RDR
=
DR\,$.
\par
Consequently:\ \ \ If $\nu (RDR) =\nu (\frak m)$  and  $T=0$  then $R$
is
regular.
\endproclaim
This was also stated in theorem 4.1  of \cite{Yo}, but there the proof
of proposition~3.3, which is used in the proof,  is wrong.
\head
3. Conditions on the embedding dimension, the index of stability,
and the multiplicity
\endhead
\subhead 3.1. Low embedding dimension
\endsubhead
In the case of  $n={\edim}\,R=2$ the ring $R$ represents  a plane curve
singularity  and therefore $R$ is a complete intersection. Then the
conjecture is true by theorem 1. But also in the cases $n=3$ and  $n=4$
one has the following results by Herzog \cite{He2, Satz 3.2 and Satz
3.3.},
that are obtained using properties of the Koszul complex for which
we refer the reader to \cite{He2}.
\proclaim{Theorem 7}
\newline
a) If $n\le 3$  and  $T=0$ then $R$ is regular.
b) If $n=4$ , $R$ is Gorenstein and $T=0$ then $R$ is regular.
\endproclaim
\subhead{3.2. Low index of stability}
\endsubhead
The main tool in this section is a reduction to the case of
${\dim}\,R=0$ ,
which was first used by Scheja in \cite{Sch}. The proofs of the
following
results are very technical, and so we will mostly contend ourselves with
references to the literature.
\newline
First we need a technical lemma which generalizes the well known formula
for
$\mu (RdR)\:$
\proclaim{Lemma 3} (\cite{G\"u2, Lemma 1a)\,} ):
Let $B=k\ldoublebracket X_1,\dots ,X_n\rdoublebracket /\frak a$  such
that
\newline
$\frak a\subseteq (X_1,\dots,X_n)^t$ for a $t\in\openN$, \ $\frak n$ the
maximal ideal of $B$ , and $\partial$  the universally finite derivation
of
$B$  over $k$.
\newline
Then for $r=t$  and $r=t+1$  we have
$$l_B (B\partial {\frak n}^r/{\frak n}^r\cdot B\partial
B)={\mu}_B({\frak n}^r) .$$
\endproclaim
This lemma together with lemma 2  is the main tool for proving
\proclaim{Proposition 4} (\cite{G\"u2, Satz 2}):
Let $\bar R$ with maximal ideal $\bar\frak m$  be an analytic
\newline
$k$-algebra
with ${\dim}\,\bar R = 0$, and $\bar d$ the universally finite
derivation of
$\bar R$ over $k$ . Assume that
$\bar R = k\ldoublebracket X_1,\dots ,X_n\rdoublebracket /\frak  a$
with
$\frak a \subseteq (X_1,\dots X_n)^r$ for some $r\in \openN$. Then
$$l_{\bar R}(\bar R \bar d \bar R)\ge (r-1)\cdot\binom{n+r-1}{r}
+\mu_{\bar R}({\bar\frak m}^r)+\mu_{\bar R}({\bar\frak m}^{r+1})\,.$$
In particular one has always (with $n = {\edim}\,\bar R$)
$$l_{\bar R}(\bar R \bar d \bar R)\ge \frac{1}{2}\cdot n\cdot (n+1)
+\mu_{\bar R}({\bar\frak m}^2)+\mu_{\bar R}({\bar\frak m}^3)\,.$$
\endproclaim
The following easy lemma enables us to apply the preceding results
to the torsion problem for one-dimensional analytic $k$-algebras:
\proclaim{Lemma 4} (\cite{Ul1, Hilfssatz 8}):
Let $R$ be a one-dimensional Cohen-Macaulay ring, $y\in\frak m$ a non
zero divisor for $R$\,, $M$  a finitely generated $R$-module with rank
$r$\, , and $T_y\:=\{ z | z\in M, y\cdot z=0\}.$
Then $l(M/y\cdot M)=r\cdot l(R/y\cdot R)+l(T_y)$.
\endproclaim
One can now show the following theorem, which is a generalization
\newline
of \cite{Ul2, Satz~7}:
\proclaim{Theorem 8} (\cite{G\"u2, Satz 4}):
Let $n\ge 3$ , and let $T_y\:=\{z|z\in T,y\cdot z=0\}$ .  If there
exists a non zero divisor $y\in\frak m$  such that ${\frak m}^4\subseteq
R\cdot
y$,
then:
$$
\alignat 2
\text{(1) If } {\edim}(R/R\cdot y) & =n-1,\ \text{then}  &\qquad
l(T_y)&=\frac{1}{2}\cdot (n-2)\cdot (n-1) \\
\text{(2) If } {\edim}(R/R\cdot y) & = n ,\ \text{then}  &\qquad
l(T_y)&=\frac{1}{2}\cdot (n-2)\cdot (n+1)\,.
\endalignat
$$
Consequently: If there exists a non zero divisor $y\in\frak m$ with
${\frak m}^4\subseteq R\cdot y$  then:
\newline
If $T=0$ then $R$ is regular.
\endproclaim
\remark{Remark 4} (\cite{G\"u2}):
\par
1) If one weakens the hypothesis of theorem 8 to ${\frak m}^5 \subseteq
R\cdot
y\,,$
one can still show $l(T_y)\ge \frac{1}{2}\cdot (n-2)\cdot (n-1)-r(R)\,,$
where $r(R)$  denotes the type of  $R$ .  So, if $R$ is Gorenstein and
$n\ge 4$ one has again $T\ne 0\, .$
\par
2) The condition ${\frak m}^{t+1}\subseteq R\cdot y$  for some $t$ is
satisfied for instance if ${\frak m}^t$ is stable in the sense of
\cite{H-W1,definition 1.2.}
(See \cite{H-W1, remark 1.5}).
\newline
If $m(R)\le {\edim}\,R + 1$  then even $ {\frak m}^3\subseteq R\cdot y$
(\cite{Ul1, Bemerkung 9 b)}): Take for $y$ a superficial element of
degree
one. Then $l_R (R/R\cdot y) = m(R)$  and so
\newline
$l(({\frak m}^2+R\cdot
y)/R\cdot y)
=l(R/R\cdot y)-l(R/({\frak m}^2+R\cdot y))=m(R) - {\edim}\, R\le 1\,.$
It follows $\frak m\cdot ({\frak m}^2 +R\cdot y)\subseteq R\cdot y $
and therefore ${\frak m}^3\subseteq R\cdot y\,.$
\endremark
This is a first example of a condition between the multiplicity
and the embedding dimension of $R$ , which will be generalized in
theorem 9.
\newline
Obviously the condition ${\frak m}^{t+1}\subseteq R\cdot y$  plays
an important role. If $y$ is a superficial element of degree $1$ then
for all large $t\in \openN$  we have $y\cdot{\frak m}^t ={\frak
m}^{t+1}\, .$
\definition{Definition 4} The minimal $t\in \openN$ such that there is a
superficial element (of degree $1$ )  $y$ with $y\cdot{\frak m}^t
={\frak
m}^{t+1}$
is called the {\it index of stability}  of $R$  and is denoted by
$t(R)\,.$
\enddefinition
\subhead
3.3. Relatively low multiplicity
\endsubhead
\proclaim{Lemma 5} (\cite{G\"u2, Lemma 3}):
Let $t(R)\ge 2$ . Then for every superficial element $x$ of degree one
$l((Rdx+x\cdot RdR)/x\cdot RdR)\ge 2.$
\endproclaim
Using lemma 4 and 5 together with proposition 4 one can now
derive:
\proclaim{Theorem 9} (\cite{G\"u2, Satz 5' and 5}):
\newline
Let $x$ be a superficial element of degree one, $\bar R\:= R/R\cdot x$ ,
and
$\bar\frak m$  the maximal ideal of $\bar R$ .
\newline
If $m(R)\le {1\over 2}\cdot n\cdot (n-1)+\mu_{\bar R}({\bar\frak m}^2)+
\mu_{\bar R}({\bar\frak m}^3)+1$  and $T=0$ then $R$ is regular.
\newline
More generally: Let $X\in k\ldoublebracket X_1,\dots,X_n\rdoublebracket$
be
a representative for $x$ , and assume that
$I\subseteq (X_1,\dots,X_n)^r +(X)$
for an  $r\in\openN$, $r\ge 2$ .
Then:
\newline
If $m(R)\le (r-1)\cdot\binom{n+r-2}{r}+\mu_{\bar R}({\bar\frak m}^r)
+\mu_{\bar R}({\bar\frak m}^{r+1})+1$  and $T=0$ then $R$ is regular.
\endproclaim
\demo{Proof}: Assume $R$  not regular ($n\ge 2$).  By lemma 4  we have
\newline
$l(T_x)=l(RdR/x\cdot RdR)-l(\bar R)=l(\bar R\bar d\bar R)+l(Rdx+x\cdot
RdR/x\cdot RdR)-l(\bar R).$
Since $x$ is superficial we have $l(\bar R)=m(R)$. By hypothesis we can
write
$\bar R=R/R\cdot x=k\ldoublebracket Z_1,\dots,Z_{n-1}\rdoublebracket
/\frak b$
with $\frak b\subseteq (Z_1,\cdots,Z_{n-1})^r$ and so by proposition 4
$l(\bar R\bar d\bar R)\ge (r-1)\cdot\binom{n+r-2}{r}+l({\bar\frak
m}^r)+l({\bar\frak m}^{r+1}).$
Together with lemma 5  (we may assume  $t(R)\ge 2$ because of theorem 7)
we get
$$l(T)\ge l(T_x)\ge (r-1)\cdot\binom{n+r-2}{r}+l({\bar\frak
m}^r)+l({\bar\frak
m}^{r+1})
+2-m(R),$$
and if the hypothesis of the theorem is satisfied, the right hand side
is positive.
\enddemo
\remark{Remark 5} (1) In \cite{Is, theorem 1} Isogawa gives a simple
proof
in the special case of the above theorem that $ m(R)\le {1\over 2}\cdot
n\cdot
(n-1)-1.$
\newline
(2) In \cite{Po1, Satz 3.13} Pohl shows that theorem 9  can be
applied for instance in the following situation: Let (with the notation
in
the proof of theorem 9)
\newline
$a\:= min\{ i | (Z_1,\dots,Z_{n-1})^i\subseteq\frak b\}$
and $c\:= max\{ j | (Z_1,\dots,Z_{n-1})^j\subseteq\frak b\}.$
\newline
If $a-c\le
2$ or
$R$ a Gorenstein singularity and $a-c\le 3$ then the hypothesis of
theorem 9  is satisfied with  $r=c.$  In particular the first condition
is
satisfied
if $R$ has {\it maximal Hilbert function}
(i.e. $\mu_R({\frak m}^i)=min\{\binom{n+i-1}{i}, m(R)\}$  for all
$i\in\openN$).
\endremark
Combining theorem 9  with  \cite{He2, Satz 3.2 and Satz 3.3} one gets:
\proclaim{Theorem 10} (\cite{G\"u2, Satz 6}): Let  R  be a
domain.
\newline
a)\ \ If $m(R)\le 9$ and $T=0$ then $R$  is regular.
\newline
b)\ \ If $m(R)\le 13$, $R$ Gorenstein, and $T=0$ then $R$ is regular.
\endproclaim
\demo{Proof} Let $x$ be a superficial element of degree one. Assume
${\frak m}^3\nsubseteq R\cdot x$  in view of theorem 8. Using the same
notation as in proposition 4 we then have $\mu({\bar\frak m}^2)\ge 1$
and
$\mu({\bar\frak m}^3)\ge 1$ .  Now either
$m(R)\le {1\over 2}\cdot n\cdot(n-1)+3\le {1\over 2}\cdot
n\cdot(n-1)+\mu({\bar\frak m}^2)
+\mu({\bar\frak m}^3)+1$, and then $R$ is regular by theorem 8, or
$m(R) > {1\over 2}\cdot n\cdot(n-1)+3 $.
\newline
On the other hand we have by hypothesis
\newline
a)\ \ $m(R)\le 9$ and therefore $n\le 3$ .\cite{He2, Satz 3.3}  gives
our
result.
\newline
b)\ \ $m(R)\le 13$ and therefore $n\le 4$ .
\newline
Now the result follows with
\cite{He2, Satz 3.2}.
\enddemo
\head
4. Conditions on the linkage class
\endhead
In this section let $k$ be algebraically closed.
\newline
Recall the definition of linkage:  Two perfect ideals $I$ and $J$
of the same grade in a Gorenstein local ring are said to be
{\it linked} (or  {\it 1-linked}), if there exists an ideal $G$
generated by
a regular sequence,\ \ $G\subset I$, $G\subset J$ , such that $I=G:J$
and $J=G:I$.
\newline
The analytic $k$-algebra~$R$ is said to be {\it linked}  (or  {\it
1-linked})
to the analytic
\newline
$k$-algebra~$R'$, if $R=k\ldoublebracket
X_1,\dots,X_n\rdoublebracket /I $
and $R'=k\ldoublebracket X_1,\dots,X_n\rdoublebracket /J$,
\newline
with $I$ linked to
$J$ .
\newline
$R$ is said to be {\it in the same linkage class} as $R'$ (or
{\it
t-linked
to  $R'$}), if there exists a sequence  $R=R_0,\dots,R_t=R'$  such
that $R_i$ is linked to $R_{i+1}$  for $i=0,\dots,t-1$.
\newline
One defines an invariant
$$\sigma(R)\:= l({\frak K}^{-1}/R)-l(R/\frak K)$$
for any canonical ideal $\frak K\subseteq R$ of $R$ . One has
$\sigma(R)\le l(S/R)=:\delta(R)$ and $\sigma(R)< \delta(R)$  if $R$ is
singular
(\cite{J\"a, Satz 1}, \cite{H-W2, pp.337/338}). If $R$ is Gorenstein
(e.g.
complete intersection) then $\sigma(R)=0$ .
\proclaim{Theorem 11} (\cite{H-W2, theorem}):
Let $R$ be 1-linked to $R'$ .Then with the notation of  1.1.:
$$l({\Coker}\, c_R)+l(\tau(I/I^2))-l(T) =\sigma(R').$$
\endproclaim
For the proof see \cite{H-W2}. Here we will look at some of the
consequences:
\newline
It follows immediately that both sides of the equation are
even-linkage invariants (that is invariants for $t$-linkage with $t$
even).
\newline
By theorem 11 one has
$$
\align
l(T)& =l({\Coker}\, c_R)+l(\tau(I/I^2)-\sigma(R') \\
    & =l(SDS/RDR)+l(\tau(I/I^2)+l(S/R)-\sigma(R').
\endalign
$$
Now, if $\sigma(R')=\sigma(R)$ and $R$ is not regular then
$l(S/R)-\sigma(R')>0$ and hence $T\ne 0$.
\newline
If $R$ and therefore also $R'\:= R_1$ is in the linkage class of a
complete intersection, it is easily seen that $R$ is also {\it evenly}
linked to a complete intersection (\cite{H-W2, p.336}). Since $\sigma$
is an
even-linkage invariant one has $\sigma(R')=0$  and therefore by the
formula of theorem 11:
\proclaim{Theorem 12} (\cite{H-W2, corollary 3}, \cite{H-W3, theorem
4.5}):
If $R$ is in the linkage class of a complete intersection then
$$
\align
l(T)&=l({\Coker}\, c_R)+l(\tau(I/I^2)) \\
    &=l(S/R)+l(SDS/RDR)+l(\tau(I/I^2)).
\endalign
$$
Consequently:\ \ \ If $T=0$ then $R$ is regular.
\endproclaim
There is a second class of rings for which this method gives a result:
\definition{Definition 5} (\cite{H-W2}):
Let $R_0$ be a reduced complete analytic $k$-algebra. $R$ is called a
{\it small extension} of $R_0$ if there exists a non zero divisor $x$
of the integral closure of $R_0$ such that $R=R_0[x]$ and $x^2\in R_0$ .
\enddefinition
Now let $R$ be a small extension of a complete intersection $R_0$ .
It is easily seen that $R$ is linked to itself, and so any $R'$ in
the linkage class of $R$ is evenly linked to $R$ . So
$\sigma(R')=\sigma(R)$
and the formula of theorem 11 yields:
\proclaim{Theorem 13} (\cite{H-W2, corollary 4}): If $R$ is in the
linkage
class of a small extension of a complete intersection then
$$
\align
l(T)&=l({\Coker}\, c_R)+l(\tau(I/I^2))-\sigma(R) \\
    &>l(SDS/RDR)+l(\tau(I/I^2)) \text{ if } R \text{ is not regular.}
\endalign
$$
Consequently:\ \ \ If $T=0$ then $R$  is regular.
\endproclaim
Quite generally one can show that of two linked singular analytic
$k$-algebras at least one of them must have non trivial torsion of
the differential module:
\proclaim{Theorem 14} (\cite{H-W2, corollary 5 and theorem on page
335}):
\newline
Let $\tau(\Omega_{R/k})=0$, $R'$  singular and 1-linked to $R$. Then
$\tau(\Omega_{R'/k})\ne 0$ .
\endproclaim
In order to show this we first need:
\proclaim{Lemma 6} (\cite{H-W2, corollary 2a}): Let $R$  be 1-linked to
a singular $R'$ and let $\delta(R)\ge\sigma(R')$ then
$\tau(\Omega_{R/k})\ne 0$ .
\endproclaim
\demo{Proof} By theorem 11 one has
$l(T)=l(SDS/RDR)+l(\tau(I/I^2))+\delta(R)-
\sigma(R')
\ge l(SDS/RDR)+l(\tau(I/I^2))$. It follows that $T\ne 0$ if $SDS\ne
RDR$.
Assume now that $SDS=RDR$, then also $S\cdot RDR=RDR$. We can choose
$s$ as in 1.2 . Then by proposition 1 $l(T)=l(R\delta R)-l(S/R)$ and
therefore $R$ has minimal torsion. By proposition 2  we then have $T\ne
0$ .
\enddemo
\demo{Proof of theorem 14} Let $R'$ be singular and 1-linked to $R$ .
Since $\tau(\Omega_{R/k})=T=0$ we must have $\sigma(R')>\sigma(R)$
because of
lemma 6. So $ \delta(R')\ge\sigma(R')>\delta(R)\ge\sigma(R)$.
Applying again
lemma 6 with $R$ and $R'$ interchanged yields
$\tau(\Omega_{R'/k})\ne 0$.
\enddemo
\head
 5. Smoothability conditions
\endhead
\definition{Definition 6} $R$ is called {\it smoothable}, if $R\cong
P/J$,
where $P$ is a normal analytic $k$-algebra and $J$ is generated by a
regular
sequence of $P$ .
\enddefinition
In the terminology of deformation theory that means that $R$ can
be deformed into a complete intersection.
\newline
Bassein \cite{Ba} shows in the complex analytic case, generalizing a
result of Pinkham \cite{Pi} for Gorenstein singularities:
\proclaim{Theorem 15} (\cite{Ba, theorem 2.4}):
\newline
If $R$ is smoothable then $l(T)=\delta (R)+{\dim}(SDS/RDR)$.
Consequently:\ \ \ If $T= 0$ then $R$ is regular.
\endproclaim
Buchweitz and Greuel (\cite{B-G, 6.1}) weakened the condition
``smoothable'' to the condition, that the degree of singularity in
the general fiber is at most ${1\over 2}\cdot\delta (R)$.
In his Dissertation Koch proves a version for an arbitrary ground
field $k$. He generalizes the above results so that in the complex
analytic case he only requires that $R$  can be deformed into an
almost complete intersection:
\proclaim{Theorem 16} (\cite{Ko, Korollar 3 to Satz 10}):
\newline
Let $P$ be excellent and an almost complete intersection in
codimension 1, and
\newline
$R=P/J$  with $J$  generated by a regular sequence
of $S$.
\newline
Then:\ \ $l(T)\ge \delta(R)+l(SDS/RDR)$.
\newline
Consequently:\ \ \ If $T= 0$ then $R$ is regular.
\endproclaim
For more results concerning smoothability see \cite{H-W3}, where
formulas involving the higher (analytic) derived functors
$T_i(R/s,R)$  and $T^i(R/s,R)$ are obtained. In particular theorem 12
is proved there by showing that $R$ is smoothable.
\head 6. Quadratic transforms
\endhead
For simplicity let us assume in this section that $R$ is a domain
and $k$ algebraically closed. Instead of looking at $T$ as the
kernel of the natural map $RdR \rightarrow SDS$, one can also take an
intermediate ring $A$  with  $R \subseteqq A \subseteqq S $  instead of
$S$ as
long as one can be sure that there is a nontrivial kernel of the induced
map of the differential modules. A good candidate for this is the first
quadratic transform
$ R_1 := R\left[ \frac{x_2}{x_1},...,\frac{x_n}{x_1} \right] $ , where
we
assume
without loss of generality that $ x_1 $  is a superficial element of
degree one
in $R$, because a short direct computation shows that whenever $ T\ne
0$  then
there is also $ {\Ker} \, \varphi \ne 0$ . The converse is trivial.
Therefore:
\proclaim
{Proposition 3} (\cite{Be3, theorem 1}):
Let $R_1$ be the first quadratic transform of $R$ and $ \varphi
:\Omega_{R/k}
\rightarrow
\Omega_{R_1/k} $ the natural homomorphism induced by the inclusion $
R\hookrightarrow R_1$.
Then $T \ne 0$ if and only if ${\Ker}\, \varphi \ne 0$ .
\endproclaim
The hope is that the behavior of the differential modules when
going from $R$ to a suitable $A$ might be easier accessible than
that when going from  $R$ to $S$ . One would think that the length
of the torsion should strictly decrease when going from a singular
$R$ to its quadratic transform $R_1$ . Indeed Bertin and Carbonne
show that this is so in the case of curve on a surface in
a point of the curve which is a simple point of the surface
(\cite{B-C}, theorem 2). If one wants only to show $T \ne 0$ it is
enough
to show $ {\Ker}\, \varphi \ne 0$  for some $A$ . By the same reasoning
as in
section 1  we obtain the analog of the formula in lemma 1:
$$ l({\Ker}\, \varphi )=l(AD_A A/RD_A R)+l(R\delta R)-l(A \partial
A)-l(A/R)\,,  $$
(with $\partial $  the universal derivation of $A$ over $s$ and $D_A$
the
universally finite derivation of $A$ over $k$) .
\par
Unfortunately, so far there is no method known to achieve this
goal for $R$ . But one can at least construct a ring $\tilde R$  with
$R\subseteq \tilde R\subseteq R_1$  which has the same multiplicity and
the same quadratic transforms as $R$ , for which by the above formula
one obtains $ l({\Ker}\, \tilde {\varphi}) \ge {1\over 2}\cdot m\cdot
(m-1)$ ,
with $\tilde{\varphi}$ denoting the natural homomorphism
$\Omega_{\tilde R /k} \rightarrow \Omega_{R_1/k}$.
More precisely:
\proclaim{Theorem 17} (\cite{Be3, corollary to theorem 2}):
Let $\tilde R = s+{\frak m}\cdot R_1$ . Then $l(\tau(\Omega_{\tilde R
/k}))\ge
\frac{1}{2}\cdot m\cdot (m-1)$ . Moreover $R$ has the same multiplicity
$m$
and the same quadratic transforms as $R$ .
\par
Consequently:\ \ If $\tau (\Omega_{\tilde R /k})= 0$  then $\tilde R$
(and
therefore
also $R$) is regular.
\endproclaim
\demo{Proof} By definition $ s = k\ldoublebracket x\rdoublebracket$ and
$x$ is
an element of minimal
value in  $R$ . We have $ {\frak m}\,\cdot R_1 = x \cdot R_1 $ ,
therefore  $x$
has minimal
value in $ \tilde R$  as well and hence $ m(\tilde R) = m(x) = m(R)\,.$
For
the maximal ideal $ \tilde{ \frak m}{\,}$  of $\tilde R$  one has $
\tilde{
\frak m}{\,}^2 =
x \cdot \frak m $ , therefore $ \tilde R $ has maximal embedding
dimension
 $ {\edim}(\tilde R) = m(\tilde R) = m $. Theorem 8  then yields
$ l(\tau (\Omega_{\tilde R /k}) \ge \frac{1}{2} m (m-1)$ . A more
detailed
 analysis in \cite{Be3} shows that this inequality already holds for
 $ {\Ker}\, \tilde{ \varphi} $.
\enddemo
\remark{Remark 6} The ring  $R$  is the  glueing  of $ \frak m\cdot
R_1$  over
$\frak m $
\ in the sense of Tamone \cite{Ta2}.
\endremark
\head 7. Equisingularity
\endhead
We conclude this paper with some remarks on the torsion and the
equisingularity class of $R$. Let us always assume that $R$ is
irreducible
and $k$ algebraically closed.
\par
For plane curves there are many equivalent definitions of
equisingularity.
For a report on the fundamental papers of Zariski \cite{Za2}
and the definitions see \cite{Ca}, Chapter III-V.
\par
The three main equisingularity classes are:
\newline
(ES1):\ $R$ and $R'$ are equisingular if they have the same
saturation.
\par\noindent
(ES2):\ $R$ and $R'$ are equisingular if they have the same multiplicity
sequence with respect to quadratic transforms.
\par\noindent
(ES3):\ $R$ and  $R'$ are equisingular if they have the same value
semigroup.
\par
While for plane curves all three concepts are equivalent (and also
equivalent to other definitions via characteristic exponents or
characteristic pairs), already for curves in 3-space
over the complex numbers there is no relation between the three
definitions (see \cite{Ca}, examples 5.3.4 and remark 5.4.4.).
But, nevertheless:
\proclaim
{Theorem 18} (\cite{Po1, Satz 2.6}):
For every singular $R$  there is always an $R'$
\newline
equisingular to $R$ with $ \tau(\Omega_{R'/k}) \ne 0$,
no matter which of the three definitions one takes.
\endproclaim
\demo{Proof}:
\par
(ES1): There is a generic plane projection $R'$ of $R$
equisingular to  $R$ , and $\tau(\Omega_{R'/k})\ne  0$  because of
$ {\edim} R' = 2$.
\par
(ES2): Take  for $R'$  the ring $\tilde R$  defined in section 6. By
theorem 17  it is (ES2)-equisingular to $R$  and $\tau(\Omega_{R'/k})
\ne 0 $ .
($R'$
has not only the same multiplicity sequence as $R$ but even the same
quadratic
transforms.)
\par
(ES3): Take for $R'=k\ldoublebracket H\rdoublebracket$ the analytic
semigroup
ring for the value semigroup
$H$  of  $R$ . By definition $R'$  is (ES3)-equisingular to  $R$  and by
remark~3
and theorem~5 we have $\tau(\Omega_{R'/k}) \ne 0 $.
\enddemo
Unfortunately the length of the torsion is not an invariant of the
equisingularity class of the curve. For instance Zariski shows in
\cite {Za1} that for a plane curves $ l(T)=2\delta (R)$  if and only if
it
can be represented by an equation $ Y^p - X^q=0$ with $ (p,q)=1$ , but
of
course there are many other plane curves $R'$ with just one
characteristic
pair$(p,q)$  (and therefore in the same equisingularity class as $R$ )
for which ${l(T) < 2\delta (R') = 2\delta  (R) }$.

\parindent=0pt

\Refs
\widestnumber\key{H-W1}

\ref\key Az
\by  A. Azevedo
\book The Jacobian ideal of a plane algebroid curve
\bookinfo Ph.D. Thesis, Purdue Univ.
\yr 1967
\endref

\ref\key Ba
\by R. Bassein
\paper On smoothable curve singularities: local methods
\jour Math. Ann
\vol 230 \pages 273--277 \yr 1977
\endref

\ref\key B-C
\by  J. Bertin et Ph. Carbonne
\paper Sur le sous-module de torsion du module des diff\'erentielles
\jour C. R. Acad. Sc. Paris, Ser. A
\vol 277 \pages 797--800 \yr 1973
\endref

\ref\key Be1
\by  R. Berger
\book \"Uber verschiedene Differentenbegriffe
\bookinfo Sitzungsber. d. Heidelberger Akad. d. Wiss., Math-naturw. Kl.,
1.
Abh.
\publ Springer--Verlag
\yr 1960
\endref

\ref\key Be2
\bysame % R. Berger
\paper Differentialmoduln eindimensionaler lokaler Ringe
\jour Math.Z. \vol 81 \pages 326--354 \yr 1963
\endref

\ref\key Be3
\bysame % R. Berger
\paper  On the torsion of the differential module of a curve singularity
\jour Arch. Math. \vol 50 \pages 526--533 \yr 1988
\endref

\ref\key B-G
\by R. Buchweitz and G. Greuel
\paper The Milnor number and deformations of complex curve singularities
\jour Inventiones Math. \vol 52 \pages 241--281 \yr 1980
\endref

\ref\key Ca
\by A. Campillo
\book Algebroid Curves in Positive Characteristic
\bookinfo Springer Lect. Notes in Math.\vol 813 \yr 1980
\publ Springer
\publaddr Berlin-Heidelberg-New York
\endref

\ref\key G\"u1
\by K. G\"uttes
\paper Einige Untersuchungen zum Torsionsproblem bei
Kurvensingularit\"aten
\book Dissertation, Saarbr\"ucken \yr 1988
\endref

\ref\key G\"u2
\by K. G\"uttes
\paper Zum Torsionsproblem bei Kurvensingularit\"aten
\jour Arch. Math. \vol 54 \pages 499--510 \yr 1990
\endref

\ref\key He1
\by J. Herzog
\paper Eindimensionale fast-vollst\"andige Durchschnitte sind nicht
starr
\jour manus. math \vol 30 \pages 1--20 \yr 1979
\endref

\ref\key He2
\by J. Herzog
\paper Ein Cohen-Macaulay Kriterium mit Anwendungen auf den
Konormalenmodul
und den Differentialmodul
\jour Math. Z.\vol 163 \pages 149--162 \yr 1978
\endref

\ref\key H-K
\by J. Herzog and E. Kunz, Editors
\book Der kanonische Modul eines Cohen-Macaulay-Rings
\bookinfo Springer Lect. Notes in Math.\vol 238 \yr 1971
\publ Springer
\publaddr Berlin-Heidelberg-New York
\endref

\ref\key H-W1
\by J. Herzog and R. Waldi
\paper A note on the Hilbertfunction of a one-dimensional
Cohen-\-Macaulay
ring
\jour Manus. Math. \vol 16 \pages 251--260 \yr 1975
\endref

\ref\key H-W2
\by J. Herzog and R. Waldi
\paper Differentials of linked curve singularities
\jour Arch. Math. \vol 42 \pages 335--343 \yr 1984
\endref

\ref\key H-W3
\by J. Herzog and R. Waldi
\paper Cotangent functors of curve singularities
\jour manuscripta math. \vol 55 \pages 307--341 \yr 1986
\endref

\ref\key H\"u
\by R. H\"ubl
\paper A note on the torsion of differential forms
\jour Arch. Math. \vol 54 \pages 142--145 \yr 1990
\endref

\ref\key Is
\by S. Isogawa
\paper On Berger's conjecture about one dimensional local rings
\jour Arch. Math. \vol 57 \pages 432--437 \yr 1991
\endref

\ref\key J\"a
\by J. J\"ager
\paper L\"angenberechnung und kanonische Ideale in eindimensionalen
Ringen
\jour Arch. Math. \vol 24 \pages 504--512 \yr 1977
\endref

\ref\key Ko
\by J. Koch
\book \"Uber die Torsion des Differentialmoduls von
Kurvensingularit\"aten
\bookinfo  Dissertation Regensburg, Regensburger Mathematische
Schriften \vol 5 \yr 1983
\publ Fakult\"at f. Math Univ. Regensburg
\endref

\ref\key Ku1
\by E. Kunz
\paper The conormal module of an almost complete intersection
\jour Proc. Amer. Math Soc. \vol 73 \pages 15--21 \yr 1979
\endref

\ref\key Ku2
\by E. Kunz
\book K\"ahler Differentials
\bookinfo Advanced Lectures in Mathematics
\publ Vieweg
\newline
\publaddr Braun\-schweig/Wies\-baden \yr 1986
\endref

\ref\key K-R
\by E. Kunz and W. Ruppert
\paper Quasihomogene Singularit\"aten algebraischer Kurven
\jour Manus. Math. \vol 22 \pages 47--61 \yr 1977
\endref

\ref\key K-W
\by E. Kunz and R. Waldi
\book Regular Differential Forms
\bookinfo Contemporary Mathematics \vol 79
\publ ASM \yr 1988
\endref

\ref\key Pi
\by H. Pinkham
\paper Deformations of algebraic varieties with $G_m$ action
\jour Ast\'erisque \vol 20 \yr 1974
\endref

\ref\key Po1
\by Th. Pohl
\book  \"Uber die Torsion des Differentialmoduls von
Kurvensingularit\"aten
\bookinfo Dis\-ser\-ta\-tion Saar\-br\"ucken \yr 1989
\endref

\ref\key Po2
\bysame %Th. Pohl
\paper Torsion des Differentialmoduls von Kurvensingularit\"aten
mit maximaler Hilbertfunktion
\jour Arch. Math. \vol 52 \pages 53--60 \yr 1989
\endref

\ref\key Po3
\by Th. Pohl
\paper Differential Modules with maximal Torsion
\jour Arch. Math. \vol 57 \pages 438--445 \yr 1991
\endref

\ref\key Sch
\by G. Scheja
\book Differentialmoduln lokaler analytischer Algebren
\bookinfo Schriftenreihe Math. Inst. Univ. Fribourg
\publ Univ. Fribourg, Switzerland
\yr 1970
\endref

\ref\key St
\by U. Storch
\paper Zur L\"angenberechnung von Moduln
\jour Arch. Math. \vol 24 \pages 39--43 \yr 1973
\endref

\ref\key Ta1
\by G. Tamone
\paper Sugli incollamenti di ideali primari e la genesi di certi
singolarit\`a
\jour Analizi Funzionale e Applicazione B.U.M.I. (Supplemento) Algebra
 e Geometria Suppl. \vol 2
\pages 243--258 \yr 1980
\endref

\ref\key Ta2
\bysame % G. Tamone
\paper Blowing--up and Glueings in one--dimensional Rings
\inbook Commutative Algebra, Proceedings of the Trento Conference
\eds S. Greco and G. Valla
\bookinfo Lect. Notes in Pure and Appl. Math \vol 84
\publ Marcel Decker, Inc. \publaddr New York and Basel
\pages 321-337 \yr 1983
\endref

\ref\key Ul1
\by B. Ulrich
\book Torsion des Differentialmoduls und Kotangentenmodul
von Kurven\-singularit\"aten
\bookinfo Dissertation Saarbr\"ucken \yr 1980
\endref

\ref\key Ul2
\by B. Ulrich
\paper Torsion des Differentialmoduls und Kotangentenmodul
von Kurven\-singularit\"aten
\jour Arch. Math \vol 36 \pages 510--523 \yr 1981
\endref

\ref\key Yo
\by Y. Yoshino
\paper Torsion of the differential modules and the value semigroup
of one dimensional local rings
\jour Math. Rep. Toyama Univ. \vol 9 \pages 83--96 \yr 1986
\endref

\ref\key Za1
\by O. Zariski
\paper  Characterization of plane algebroid curves whose module of
differentials has maximum torsion
\jour Proc. Nat. Acad. Sci. \vol 56 \pages 781--786 \yr 1966
\endref

\ref\key Za2
\bysame % O. Zariski
\paper Studies in Equisingularity I
\jour Amer. J. Math. \vol 87 \pages 507--535 \yr 1965
\moreref \paper II \vol 87 \pages 972--1006 \yr 1965
\moreref \paper III \vol 90 \pages 961--1023 \yr 1965
\endref

\endRefs
\enddocument